# The Role of Kemeny's Constant in Properties of Markov Chains

Jeffrey J. Hunter[1]

*School of Computing and Mathematical Sciences, Auckland University of Technology, Auckland, New Zealand*

In a finite irreducible Markov chain with stationary probabilities $\{\pi_i\}$ and mean first passage times $m_{ij}$ (mean recurrence time when $i = j$) it was first shown, by Kemeny and Snell (1960) that $\sum_j \pi_j m_{ij}$ is a constant, $K$, (Kemeny's constant) not depending on $i$. A variety of techniques for finding expressions and bounds for $K$ are given. The main interpretation focuses on its role as the expected time to mixing in a Markov chain. Various applications are considered including perturbation results, mixing on directed graphs and its relation to the Kirchhoff index of regular graphs.

**Introduction**

In a finite $m$-state irreducible Markov chain with stationary probabilities $\{\pi_i\}$ and mean first passage times $m_{ij}$ (mean recurrence time when $i = j$) Kemeny and Snell, (1960), showed that $\sum_{j=1}^{m} \pi_j m_{ij}$ is a constant, $K$, not depending on $i$. This constant has since become known as *Kemeny's constant*. Interest in this constant is growing in that to date no reasoned argument, apart from mathematical justifications, has been put forward as a justification as to why this feature of a Markov chain should be a constant.

We first introduce some background theory that leads to various expressions for Kemeny's constant. We then follow, basically, a chronological development of the exploration of the derivation of the constant and give some interpretations. Applications to perturbed Markov chains, random walks on directed graphs and a recent relationship to the Kirchhoff index on regular graphs through the consideration of the properties of random walks via electric networks is considered. This paper is a survey paper with the inclusion of some new derivations and expressions. The link with the Kirchhoff index raises the opportunity for further research to clarify possible alternative linkages.

---

Address for correspondence: Jeffrey J Hunter, School of Computing and Mathematical Sciences, Auckland University of Technology, Private Bag 92006, Auckland 1142, New Zealand. Email: Jeffrey.hunter@aut.ac.nz

**Markov chain theory**

Let $X_n$, $(n \geq 0)$ be a finite irreducible (ergodic), discrete time Markov chain (MC). Let $S = \{1, 2, \ldots, m\}$ be its state space. Let $p_{ij} = P\{X_{n+1} = j | X_n = i\}$ be the transition probability from state $i$ to state $j$. Let $P = [p_{ij}]$ be the transition matrix of the MC. Since $P$ is a stochastic matrix, $\sum_{j=1}^{m} p_{ij} = 1$ for all $i \in S$. Let $\{p_j^{(n)}\} = \{P[X_n = j]\}$ be the probability distribution of $X_n$ at the $n$-th trial. When the MC is finite, aperiodic and irreducible, i.e. regular, a limiting distribution exists that does not depend on the initial distribution and the limiting distribution is the stationary distribution. i.e. $\{X_n\}$ has a unique stationary distribution $\{\pi_j\}, j \in S$ and $\lim_{n \to \infty} p_j^{(n)} = \pi_j$. When the MC is finite, irreducible and periodic, a limiting distribution does not exist. However it does have a unique stationary distribution. In summary, irreducible MCs $\{X_n\}$ have a unique stationary distribution $\{\pi_j\}, j \in S$, with stationary probabilities $\pi_i$ given as the solution of the stationary equations: $\pi_j = \sum_{i=1}^{m} \pi_i p_{ij}$ $(j \in S)$ with $\sum_{i=1}^{m} \pi_i = 1$. We represent the stationary probability vector as $\boldsymbol{\pi}^T = (\pi_1, \pi_2, \ldots, \pi_m)$.

We present a primer on g-inverses of $I - P$. A "one condition" g-inverse, or an "equation solving" g-inverse, of a matrix $A$ is any matrix $A^-$ such that $AA^-A = A$. These g-inverses are used extensively in solving systems of linear equations. A necessary and sufficient condition for $AXB = C$ to have a solution is that $AA^-CB^-B = C$. If this consistency condition is satisfied, the general solution is given by $X = A^-CB^- + W - A^-AWBB^-$, where $W$ is an arbitrary matrix, (Rao, (1966)). In particular, if the consistency condition $AA^-C = C$,

$$AX = C \qquad (1)$$

has a solution $X = A^-C + (I - A^-A)W$ where $W$ is an arbitrary matrix.

This theory has been applied to Markov chains. Let $P$ be the transition matrix of a finite irreducible MC with stationary probability vector $\boldsymbol{\pi}^T$. Let $\boldsymbol{t}$ and $\boldsymbol{u}$ be any vectors and let $\boldsymbol{e}^T = (1,1, \ldots, 1)$. $I - P + \boldsymbol{tu}^T$ is non-singular if and only if $\boldsymbol{\pi}^T \boldsymbol{t} \neq 0$ and $\boldsymbol{u}^T \boldsymbol{e} \neq 0$. Further, if $\boldsymbol{\pi}^T \boldsymbol{t} \neq 0$ and $\boldsymbol{u}^T \boldsymbol{e} \neq 0$, then $[I - P + \boldsymbol{tu}^T]^{-1}$ is a g-inverse of $I - P$. (Hunter, (1982)). Actually, any g-inverse of $I - P$ can be expressed as $G = [I - P + \boldsymbol{tu}^T]^{-1} + \boldsymbol{ef}^T + \boldsymbol{g\pi}^T$ where $\boldsymbol{f}, \boldsymbol{g}, \boldsymbol{t}$ and $\boldsymbol{u}$ are vectors with $\boldsymbol{\pi}^T \boldsymbol{t} \neq 0$ and $\boldsymbol{u}^T \boldsymbol{e} \neq 0$. (Hunter, (1982)). In particular there are two special g-inverses of $I - P$ that have been used extensively in the literature: $Z \equiv [I - P + \Pi]^{-1}$, $(\Pi \equiv \boldsymbol{e\pi}^T)$, the fundamental matrix of irreducible (ergodic) MCs due to Kemeny and Snell (1960) and shown to be a g-inverse of $I - P$ by Hunter (1969), and $A^\# = (I - P)^\# = Z - \Pi$, the group inverse of $I - P$ established by Meyer, (1975).

Let us now introduce first passage times in MCs. Let $T_{ij}$ be the first passage time random variable from state $i$ to state $j$, i.e. $T_{ij} = \min\{n \geq 1$ such that $X_n = j$ given that $X_0 = i\}$, so that $T_{ii}$ is the "first return" time to state $i$. The irreducibility of the MC ensures that the $T_{ij}$ are all proper random variables. Under the finite state space restriction, all the moments of $T_{ij}$ are finite. Let $m_{ij}$ be the mean first passage time from state $i$ to state $j$. i.e. $m_{ij} = E[T_{ij} | X_0 = i]$ for all $(i,j) \in S \times S$. For an irreducible finite MC with transition matrix $P$, let $M = [m_{ij}]$ be the matrix of expected first passage times from state $i$ to state $j$. From Kemeny and Snell (1960), $M$ satisfies the matrix equation

$$(I - P)M = E - PM_d, \qquad (2)$$

where $E = ee^T = [1]$, $M_d = [\delta_{ij}m_{ij}] = (\Pi_d)^{-1} \equiv D$.

Hunter (1982) showed that equation (2), which is of the form of equation (1), can be solved using any g-inverse $G$ of $I - P$ to yield the solution
$$M = [G\Pi - E(G\Pi)_d + I - G + EG_d]D. \qquad (3)$$
Hunter (2008) showed that under any of the following three equivalent conditions: (i) $Ge = ge$, $g$ a constant, (ii) $GE - E(G\Pi)_d D = 0$, or (iii) $G\Pi - E(G\Pi)_d = 0$, that
$$M = [I - G + EG_d]D. \qquad (4)$$
Equation (4) holds in two special cases, viz. $G = Z$ (when $g = 1$) and $G = A^{\#} = Z - \Pi$ (when $g = 0$).

From equation (3), if $G = [g_{ij}]$ is any generalized inverse of $I - P$,
$$m_{ij} = \left(\frac{g_{jj} - g_{ij} + \delta_{ij}}{\pi_j}\right) + (g_{i.} - g_{j.}), \quad \text{for all } i, j. \qquad (5)$$

Under the conditions of equation (4), i.e. in particular when $G = Z$ or $A^{\#}$,
$$m_{ij} = \left(\frac{g_{jj} - g_{ij} + \delta_{ij}}{\pi_j}\right), \quad \text{for all } i, j. \qquad (6)$$

Thus for $i \neq j$, $m_{ij} = (z_{jj} - z_{ij})/\pi_j = (a^{\#}_{jj} - a^{\#}_{ij})/\pi_j$ with $m_{ii} = 1/\pi_i$ (the "mean recurrence time of state $j$") where $Z = [z_{ij}]$ and $A^{\#} = [a^{\#}_{ij}]$. These results were first derived by Kemeny and Snell (1960), using $Z$, and Meyer (1975), using $A^{\#}$.

**Kemeny's constant**

We have the following key result:

For all $i \in S$, $K_i \equiv \sum_{j=1}^{m} m_{ij}\pi_j = \sum_{j=1}^{m} z_{jj} = tr(Z) = K$, a constant not depending on $i$. (7)

$K$ is known as *Kemeny's constant*. Equivalently, in matrix-vector form, the conclusion above implies $M\pi = Ke$. A simple proof can be given, using equation (4) with $G = Z$. $M\pi = [I - Z + EZ_d]D\pi = [I - Z + EZ_d]e = e - Ze + ee^T Z_d e = Ke$, where $K = e^T Z_d e = tr(Z)$.

Since its initial appearance in Kemeny and Snell (1960), result (7) has intrigued researchers. Can we give a simple proof of its constant nature? Is there an intuitively plausible reason for this constant to arise and what are its key properties? We give a chronological development of the history of the constant while providing some alternative, clearer, proofs of the constant nature of $K$.

Before this we need to clarify that there are two versions of Kemeny's constant: Since $m_{ii} = 1/\pi_i$, $K = \sum_{j=1}^{m} m_{ij}\pi_j = m_{ii}\pi_i + \sum_{j \neq i} m_{ij}\pi_j = 1 + \sum_{j \neq i} m_{ij}\pi_j$. Some authors define, by convention, that $m_{ii} = 0$ so that the expression for the mean first passage times taken as $m_{ij} = (z_{jj} - z_{ij})/\pi_j$ holds for all $i, j$. We will stay with the expression as defined in result (7) for $K$, bearing in mind that in some books and papers $K$ is replaced by $K - 1$, which we shall call the modified Kemeny's constant.

The initial appearance of result (7) was in the book by Kemeny and Snell (1960) (p. 81) as a stated theorem with no interpretation and no implications as to its relevance.

The next appearance of $K$ was in the book by Snell (1975) (pp. 289-290) in a problem where he seeks a derivation of result (7) (with $K$ replaced by $K - 1$) and states "A prize is offered for the first person to give an intuitively plausible reason for the above sum to be independent of $i$."

This is taken up in the book by Grinstead and Snell (1997) in its first edition where it is reported (see Q24, p. 470 of the second edition (2009)) "In the course of a walk with Snell along Minnehaha Avenue in Minneapolis in the fall of 1983, Peter Doyle suggested the following explanation for the constancy of Kemeny's constant. Choose a target state according to a fixed vector $\pi$. Start from state $i$ and wait until the time $T$ that the target state occurs for the first time. Let $K_i$ be the expected value of $T$. Observe that $K_i + \pi_i \cdot 1/\pi_i = \sum_j p_{ij} K_j + 1$ and hence $K_i = \sum_j p_{ij} K_j$. By the maximum principle $K_i$ is a constant. Should Peter have been given the prize?"

We shed some more light on this by giving an alternative derivation. Define $\mathbf{k} \equiv M\pi$, where $\mathbf{k}^T = (K_1, K_2, ...., K_m)$. From equation (2), $(I - P)\mathbf{k} = (I - P)M\pi = E\pi - PM_d\pi = ee^T\pi - Pe = e - e = \mathbf{0}$. i.e. $\mathbf{k} = P\mathbf{k}$. It is easy to establish, by induction, from this observation that, for all $n \geq 1$, $\mathbf{k} = P^n\mathbf{k}$,

$$\text{i.e. } K_i = \sum_{j=1}^{m} p_{ij}^{(n)} K_j, \tag{8}$$

since the $n$-step transition probabilities of the MC, $p_{ij}^{(n)}$ are given as $\left[p_{ij}^{(n)}\right] = P^n$. By definition (7), each $K_i > 0$. Let $K^* = \max_{1 \leq i \leq m} K_i \ (> 0)$. There will be at least one state $a$ such that $K_a = K^*$. Let us assume that there is at least one state $b$ such that $0 < K_b < K^*$. The irreducibility of the MC implies that for this pair of states there exists an $n$ such that $p_{ab}^{(n)} > 0$. Thus $K^* = K_a = \sum_{j=1}^{m} p_{aj}^{(n)} K_j < \sum_{j=1}^{m} p_{aj}^{(n)} K^*$, since for all $j \neq b$, $K_j \leq K^*$ while $K_b < K^*$. This implies $K^* < \left(\sum_{j=1}^{m} p_{aj}^{(n)}\right) K^* = K^*$, a contradiction, so that there is no such state $b$, and consequently $K_i = K$ for all $i$.

As an alternative proof, note that since $\mathbf{k} = P\mathbf{k}$ the irreducibility of the MC implies that $\mathbf{k}$ is a right eigenvector of $P$ corresponding to the eigenvalue $\lambda = 1$ which must be a scalar multiple of $e$, so that $K_i = K$ for all $i \in S$.

Alternatively, if the MC is regular, (finite and aperiodic), $\lim_{n \to \infty} p_{ij}^{(n)} = \pi_j$, so that taking the limit as $n \to \infty$ in equation (8), $K_i = \lim_{n \to \infty} \left(\sum_{j=1}^{m} p_{ij}^{(n)} K_j\right) = \sum_{j=1}^{m} \pi_j K_j$, a constant $(K)$ not depending on $i$.

Doyle (2003) continues the story by stating that "a prize was offered not for an intuitive interpretation for Kemeny's constant, but for an intuitively plausible reason for $K$ to be constant. So far we've seen nothing approaching an argument for why the expected time to equilibrium should be independent of the starting state $i$, so it's hard to see why Peter should have been given the prize. On the other hand, we may note that this question does not ask whether Peter should have been given the prize on the basis of this interpretation, but just whether he should have been given the prize. Did Peter offer an intuitively plausible reason for $K_i$ to be constant, somehow related to its being the expected time to equilibrium? How can we decide this question on the basis of the information we've been given? Grinstead and Snell are not in the habit of leaving us high and dry, so the most likely explanation is that there is a simple and intuitively plausible reason for the constancy of $K_i$, immediately related to the interpretation as the time to equilibrium" (which we explore in more carefully below.)

Doyle (2009) considers the above equation $K_i = \sum_j p_{ij} K_j$ in more detail, stating "By the familiar maximum principle, any function $f_i$ satisfying $\sum_j p_{ij} f_j = f_i$ must be constant: Choose $i$ to maximise $f_i$ and observe that the maximum must be attained also for any $j$ where $p_{ij} > 0$; push the max around until it is attained everywhere. So $K_i$ doesn't depend on $i$." He then states that "the application of the maximum principle shows that the only column eigenvectors having eigenvalue 1 for the matrix $P$ are the constant vectors – a fact that was stated not quite explicitly above. The formula provides a computational verification that Kemeny's constant is a constant but doesn't explain why it is constant. Kemeny felt this keenly: A prize was offered for a more conceptual proof, and awarded – rightly or wrongly – on the basis of the maximum principle argument outlined above." As a side comment: "When Laurie Snell mailed Peter Doyle the prize for Kemeny's constant, he first made the mistake of trying to send a $50 bill by mail. That first letter never arrived. You have probably heard before that you should not send cash through the mail. Let this be a further lesson!"

Hunter (2006) gave some results that lead directly to general expressions for Kemeny's constant using g-inverses. If $G = [g_{ij}]$ is any g - inverse of $I - P$, then

$$K = 1 + tr(G) - tr(G\Pi) = 1 + \sum_{j=1}^{m} (g_{jj} - g_{j.}\pi_j).$$

When $Ge = ge$, $K = 1 - g + tr(G) = 1 - g + \sum_{j=1}^{m} g_{jj}$.

In particular, $K = tr(Z) = \sum_{j=1}^{m} z_{jj}$ and $K = 1 + tr(A^{\#})$.

These last two expressions are the usual "classical results" but the expression with $Z$ was also derived indirectly by Lovasz and Winkler, (1998), in their "random target lemma". A further derivation appears in the book on "Reversible MCs and Random walks" by Aldous and Fill (1999).

Alternative expressions for Kemeny's constant can be given in terms of eigenvalues. Since $P$ is irreducible the eigenvalues of $P$, $\{\lambda_i\}$ ($i = 1, 2,..., m$) are such that $\lambda_1 = 1$, with $|\lambda_i| \leq 1$ and $\lambda_i \neq 1$ ($i = 2,..., m$). The eigenvalues of $Z = [z_{ij}] = [I - P + e\pi^T]^{-1}$ are $\lambda_i(Z) = 1$ ($i = 1$), $\frac{1}{1-\lambda_i}$ ($i = 2,..., m$). Thus, from Levene and Loizou, (2002), Hunter, (2006), Doyle, (2009):

$$K = tr(Z) = \sum_{i=1}^{m} z_{ii} = \sum_{i=1}^{m} \lambda_i(Z) = 1 + \sum_{i=2}^{m} \frac{1}{1-\lambda_i}. \qquad (9)$$

Equation (9) forms the basis of for deriving bounds on $K$. Since $P$ is irreducible if any eigenvalue appears on the unit circle $|\lambda| = 1$ it must appear as a root of unity and be associated with a periodic chain (whose periodicity cannot exceed $m$). Any complex root $\lambda = a + bi$ must be associated with its complex conjugate $\bar{\lambda} = a - bi$, with $a^2 + b^2 \leq 1$. For this pair of conjugate roots

$$\frac{1}{1-\lambda} + \frac{1}{1-\bar{\lambda}} = \frac{2-(\lambda+\bar{\lambda})}{(1-\lambda)(1-\bar{\lambda})} = \frac{2-2a}{1-(\lambda+\bar{\lambda})+\lambda\bar{\lambda}} = \frac{2-2a}{1-2a+a^2+b^2} \geq 1.$$

For any real roots $-1 \leq \lambda \leq 1$ we must have $\frac{1}{1-\lambda} \geq \frac{1}{2}$. The only possible root at $\lambda = -1$ occurs with a periodic MC with an even period. Thus taking the real roots individually and the complex roots in pairs it is easily seen that
$$K = 1 + \sum_{i=2}^{m} \frac{1}{1-\lambda_i} \geq 1 + \frac{m-1}{2} = \frac{m+1}{2}.$$
This proof was given by Hunter (2006), based on results of Styan (1964) who derived the result for the case of $\lambda_i$ real. Thus for all finite $m$-state irreducible MCs we have a lower bound for $K$ of $(m+1)/2$.

If the MC is reversible (in which case all the $\lambda_i$ are real) and regular (aperiodic) then Levene and Loizou (2002) showed that $\frac{m-1}{2} \leq \sum_{i=2}^{m} \frac{1}{1-\lambda_i} \leq \frac{m-1}{1-\lambda_2}$.

Note also that $K = 1 + \sum_{i=2}^{m} \frac{1}{1-\lambda_i} = m + \sum_{i=2}^{m} \frac{\lambda_i}{1-\lambda_i}$.

The above bounds were improved by Palocois and Remon (2010) under the assumption that the MC is irreducible and reversible so that $1 = \lambda_1 > \lambda_2 \geq ... \geq \lambda_m > -1$. Apply the method of Lagrange multipliers to the function $f(x_2,...,x_m) = \sum_{i=2}^{m} \frac{x_i}{1-x_i}$, subject to $1 + x_2 + ... + x_m = 0$ on the domain $1 > x_2 \geq ... \geq x_m > -1$. This implies that the minimum of $f(x_1, x_2,...,x_m)$ is attained at $x_2 = .. = x_m = \frac{-1}{m-1}$, leading to the bounds
$$\frac{(m-1)^2}{m} \leq \sum_{i=2}^{m} \frac{1}{1-\lambda_i} \leq \frac{m-1}{1-\lambda_2}.$$

An alternative representation of $K$ has been given by Catral, Kirkland, Neumann, Sze, (2010), viz., $K = tr(A_j^{-1}) - \frac{A_{jj}^{\#}}{\pi_j} + 1$, where $A_j^{-1}$ is $(m-1) \times (m-1)$ principal submatrix of $A = I - P$ obtained by deleting $j$-th row and column. The proof is based upon expressing $A^{\#} = [a_{ij}^{\#}]$ in terms of $A_n^{-1}$ and $\pi^T$. Without loss of generality, take $j = m$. $m_{ij}\pi_j = a_{jj}^{\#} - a_{ij}^{\#}$ and the result follows from the result of Meyer (1973), that if $B$ is the leading $(m-1) \times (m-1)$ principal submatrix of $A^{\#}$, then $B = A_n^{-1} + \beta W - A_n^{-1}W - WA_n^{-1}$, where $\beta = \boldsymbol{u}^T A_n^{-1} \boldsymbol{e}$, $W = \boldsymbol{e}\boldsymbol{u}^T$ and $\boldsymbol{\pi}^T = (\boldsymbol{u}^T, \pi_n)$.

**Mixing Times in Markov chains**

For all irreducible MCs (including periodic chains), if for some $k \geq 0$, $p_j^{(k)} = P[X_k = j]$ $= \pi_j$ for all $j \in S$, then $p_j^{(n)} = P[X_n = j] = \pi_j$ for all $n \geq k$ and all $j \in S$. How many trials do we need to take so that $P[X_n = j] = \pi_j$ for all $j \in S$? If this is not possible, can we get arbitrarily close to $\pi_j$? Alternatively, can we choose a trial $n$ where $P[X_n = j] = \pi_j$? This is the motivation to the concept of "mixing times in MCs".

Let $Y$ be a random variable whose probability distribution is the stationary distribution $\{\pi_j\}$. We say that the MC $\{X_n\}$ achieves "mixing" at time $T = k$, when $X_k = Y$ for the smallest such $k \geq 1$ and call $T$ the "time to mixing" in the MC. Thus, we first sample from the stationary distribution $\{\pi_j\}$ to determine a value of the random variable $Y$, say $Y = j$. Now observe the MC, starting at a given state $i$. We achieve "mixing" at time $T = n$ when $X_n = j$ for the first such $n \geq 1$.

The finite state space and irreducibility of the MC implies that $T$ is finite (a.s) with finite moments (see Hunter (2006). Consider the "expected time to mixing", starting at state $i$, (assuming that mixing cannot occur at the first trial). Conditional upon $X_0 = i$, $E[T] = E_Y(E[T \mid Y]) = \sum_{j=1}^{m} E[T \mid Y = j]P[Y = j] = \sum_{j=1}^{m} E[T_{ij} \mid X_0 = i]\pi_j = \sum_{j=1}^{m} m_{ij}\pi_j = K_i = K$.

Thus the expected time to mixing, starting in any state, is $K$. (Hunter, (2006)). This formalises the "target state" or "time to equilibrium" concepts of Doyle (2003). It does not however provide a reason as to why this should be a constant.

This is all the more intriguing since if we start a finite irreducible MC in equilibrium with a state $i$ chosen according to the stationary distribution $\{\pi_i\}$ and consider the time $T_j^*$ to a fixed state $j$, this time has expectation $\sum_{i=1}^{m} \pi_i m_{ij}$. From Corollary 2.5.7 of Hunter, (2008), it can be shown that if $G = [g_{ij}]$ is any g - inverse of $I - P$, then $E(T_j^*) = \sum_{i=1}^{m} \pi_i m_{ij} = 1 + \sum_{i=1}^{m} \pi_i g_{i \cdot} - g_{j \cdot} + \left(g_{jj} - \sum_{i=1}^{m} \pi_i g_{ij}\right)/\pi_j$.

When $G\mathbf{e} = g\mathbf{e}$, $E(T_j^*) = \sum_{i=1}^{m} \pi_i m_{ij} = 1 + \left(g_{jj} - \sum_{i=1}^{m} \pi_i g_{ij}\right)/\pi_j$.

In particular, $E(T_j^*) = z_{jj}/\pi_j = 1 + a_{jj}^{\#}/\pi_j$, an expectation that is, in general, not constant!

We have defined "mixing" to occur at time $T = k$, when $X_k = Y$ for the smallest such $k \geq 1$ when the initial starting state is $i$. Suppose however we allow mixing to be possible when $k = 0$ when $i = j$. i.e. we permit "mixing" to occur at time $T = 0$, when state $i$ is the "hitting" state (rather than "returned" state). The expected time to mixing in this situation would be $\sum_{j \neq i} m_{ij}\pi_j = K - 1$, since $m_{ii}\pi_i = 1$. Thus the modified Kemeny constant has the interpretation as the expected time to "hitting" an equilibrium state. Hunter (2013) considers the distribution of the time to mixing and the time to hitting in finite irreducible MCs.

Levene and Loizou, (2002), give the following "random surfer" interpretation to Kemeny's constant. They rewrite $K$ as $\sum_{i=1}^{m} \pi_i \sum_{j=1}^{m} \pi_j m_{ij} = \sum_{i=1}^{m} \pi_i M_i$ where $M_i = \sum_{j=1}^{m} \pi_j m_{ij}$. $M_i$ gives "the mean first passage time from state $i$ when the destination state is unknown". $K = \sum_{i=1}^{m} \pi_i M_i$ can be interpreted as "the mean first passage time from an unknown starting state to an unknown destination state". They give the following interpretation: "Imagine therefore a random surfer who is following links according to the transition probabilities. At some stage our random surfer is "lost" and does not know the state he is at and where he is heading for. In this context Kemeny's constant can be interpreted as the mean number of links the random surfer follows before reaching his destination. Thus the random surfer is not "lost" anymore, he just has to follow $K$ random links and he can expect to arrive at his final destination". (As an aside, note that $M_i = K_i = K$ so that the interpretation does not lead to any reason for the constant nature of $K$.)

We conclude this section on two examples by computing the Kemeny constants. For further examples see Hunter (2006).

*Example – Two state Markov Chains*

Let $P = \begin{bmatrix} p_{11} & p_{12} \\ p_{21} & p_{22} \end{bmatrix} = \begin{bmatrix} 1-a & a \\ b & 1-b \end{bmatrix}$, ($0 \le a \le 1$, $0 \le b \le 1$). Let $d = 1 - a - b$.

The MC is irreducible if and only if $-1 \le d < 1$, in which case the MC has a unique stationary probability vector given by $\pi^T = (\pi_1, \pi_2) = \left(\frac{b}{a+b}, \frac{a}{a+b}\right) = \left(\frac{b}{1-d}, \frac{a}{1-d}\right)$.

The MC is regular if $-1 < d < 1$ and periodic, period 2 if $d = -1$ (with $a = 1$, $b = 1$)

The mean first passage time matrix is $M = \begin{bmatrix} \frac{1-d}{b} & \frac{1}{a} \\ \frac{1}{b} & \frac{1-d}{a} \end{bmatrix}$.

Kemeny's constant is given as $K = 1 + \frac{1}{a+b} = 1 + \frac{1}{1-d}$.

The minimum value of $K$ is 1.5 and occurs when the MC is periodic, period 2.
Under independent trials $b = 1 - a$, $d = 0$ and $K = 2$.
When both $a \to 0$ and $b \to 0$, so that $d \to 1$, $K$ becomes arbitrarily large.
For all two state MCs: $1.5 \le K < \infty$.

*Example – Three state Markov Chains*

Let $P = [p_{ij}] = \begin{bmatrix} 1 - p_2 - p_3 & p_2 & p_3 \\ q_1 & 1 - q_1 - q_3 & q_3 \\ r_1 & r_2 & 1 - r_1 - r_2 \end{bmatrix}$.

There are six constrained parameters with $0 < p_2 + p_3 \le 1$, $0 < q_1 + q_3 \le 1$ and $0 < r_1 + r_2 \le 1$. Let $\Delta_1 \equiv q_3 r_1 + q_1 r_2 + q_1 r_1$, $\Delta_2 \equiv r_1 p_2 + r_2 p_3 + r_2 p_2$, $\Delta_3 \equiv p_2 q_3 + p_3 q_1 + p_3 q_3$, and $\Delta \equiv \Delta_1 + \Delta_2 + \Delta_3$.

The MC is irreducible, and hence a stationary distribution exists, if and only if $\Delta_1 > 0$, $\Delta_2 > 0$, $\Delta_3 > 0$, in which case the stationary distribution given by $(\pi_1, \pi_2, \pi_3)$ $= \frac{1}{\Delta}(\Delta_1, \Delta_2, \Delta_3)$.

Let $\tau_{12} = p_3 + r_1 + r_2$, $\tau_{13} = p_2 + q_1 + q_3$, $\tau_{21} = q_3 + r_1 + r_2$, $\tau_{23} = q_1 + p_2 + p_3$, $\tau_{31} = r_2 + q_1 + q_3$, $\tau_{32} = r_1 + p_2 + p_3$. Let $\tau = p_2 + p_3 + q_1 + q_3 + r_1 + r_2$ implying $\tau = \tau_{12} + \tau_{13} = \tau_{21} + \tau_{23} = \tau_{31} + \tau_{32}$.

The mean first passage time matrix is given by

$M = \begin{bmatrix} \Delta/\Delta_1 & \tau_{12}/\Delta_2 & \tau_{13}/\Delta_3 \\ \tau_{21}/\Delta_1 & \Delta/\Delta_2 & \tau_{23}/\Delta_3 \\ \tau_{31}/\Delta_1 & \tau_{32}/\Delta_2 & \Delta/\Delta_3 \end{bmatrix}$.

Kemeny's constant is given as

$K = 1 + \frac{\tau}{\Delta} = \frac{p_2 + p_3 + q_1 + q_3 + r_1 + r_2}{q_3 r_1 + q_1 r_2 + q_1 r_1 + r_1 p_2 + r_2 p_3 + r_2 p_2 + p_2 q_3 + p_3 q_1 + p_3 q_3}$.

For all three state irreducible MCs, $K \ge 2$.
$K = 2$ is achieved in the minimal period 3 case when $p_2 = q_3 = r_1$.

In the "constant movement" case when $p_2 + p_3 = q_1 + q_3 = r_1 + r_2 = 1$,

$$K = 1 + \frac{3}{3 - q_3 r_2 - r_1 p_3 - p_2 q_1} \text{ implying } 2 \leq K \leq 2.5.$$

For the "period 2 case", when we can have transitions between {1,3} and {2}, i.e.

$$P = \begin{bmatrix} 0 & 1 & 0 \\ q_1 & 0 & q_3 \\ 0 & 1 & 0 \end{bmatrix}, (q_1 + q_3 = 1), \text{ we achieve } K = 2.5.$$

*Example – General m-state Markov Chains*
In an $m$-state MC, the minimal value of $K = (m + 1)/2$ is achieved for a periodic, period $m$, chain. For a sequence of independent trials with $m$ possible outcomes, $K = m$.
For all irreducible $m$-state MCs, $(m + 1)/2 \leq K < \infty$.

**Perturbation results**
Consider perturbing the transition matrix $P = [p_{ij}]$ associated with an ergodic, $m$-state MC, to $\overline{P} = [\overline{p_{ij}}] = P + \boldsymbol{E}$ where $\boldsymbol{E} = [\varepsilon_{ij}]$, ($\sum_{j=1}^{m} \varepsilon_{ij} = 0$). If $\overline{P}$ is also ergodic, let $\pi^T = (\pi_1, \pi_2, ..., \pi_m)$ and $\overline{\pi}^T = (\overline{\pi}_1, \overline{\pi}_2, ..., \overline{\pi}_m)$ be the associated stationary probability vectors of the two MCs. Hunter, (2006), showed that, for $m$-state ergodic MCs undergoing such a perturbation, $\|\pi^T - \overline{\pi}^T\|_1 \leq (K-1)\|\boldsymbol{E}\|_\infty$ i.e.

$$\sum_{j=1}^{m} |\pi_j^T - \overline{\pi}_j^T| \leq (K-1) \max_{1 \leq i \leq m} \sum_{k=1}^{m} |\varepsilon_{ki}|.$$

Catral, Kirkland, Neumann, and Sze, (2010) consider a variety of elementary perturbations of the above forms. Let $M = [m_{ij}]$ and $\overline{M} = [\overline{m_{ij}}]$ be the mean first passage matrices and let $K$ and $\overline{K}$ be the Kemeny constants associated with $P$ and $\overline{P}$.
For a Type 1 perturbation, when $\boldsymbol{E} = e_r h^T$, where $e_r$ is the $r$–th elementary vector, $\overline{m_{ir}} = m_{ir}$ for all $i \neq r$; $\overline{m_{ij}} \geq m_{ij} \Leftrightarrow \overline{\pi}_j \leq \pi_j$ for all $i, j \neq r$. $K \leq \overline{K}$ if and only if

$$\sum_{i \neq r} (\overline{\pi}_i - \pi_i) m_{ir} \geq 0.$$

For a Type 2 perturbation when $\boldsymbol{E} = e h^T$, $K = \overline{K}$.
They also extended their results to include:
1. Let $P$ be a symmetric stochastic, irreducible matrix $\overline{P} = P - \boldsymbol{E}$ where $\boldsymbol{E}$ is a positive semi-definite matrix with $\overline{P}$ stochastic. Then $\sum_{j=1}^{m} \overline{m_{ij}} \leq \sum_{j=1}^{m} m_{ij}$ and $\overline{K} \leq K$.
2. Let $P$ be a stochastic, irreducible matrix and suppose $0 \leq \alpha \leq 1$. Let $\overline{P} = \alpha P + (1-\alpha) e v^T$ where $v^T$ is a positive probability vector, then $\overline{K} \leq K$.

**Directed Graphs**
A directed graph, or digraph $\mathcal{G} = (\mathcal{V}, \mathcal{E})$ is a collection of vertices (or nodes) $i \in \mathcal{V} = \{1,...,m\}$ and directed edges or arcs $(i \to j) \in \mathcal{E}$. One can assign weights to each directed edge, making it a weighted digraph. An unweighted digraph has common edge weight 1. $\mathcal{G}$ can be represented by its m × m adjacency matrix $A = [a_{ij}]$ where $a_{ij} \neq 0$ is the weight on arc $(i \to j)$ and $a_{ij} = 0$ if $(i \to j) \notin \mathcal{E}$. A digraph $\mathcal{G}$ is strongly connected, or a strong digraph, if there is a path $i = i_0 \to i_1 \to ... \to i_k = j$ for any pair of nodes where each link $i_{r-1} \to i_r \in \mathcal{E}$. We focus on strong digraphs.

A random walk over a graph can be represented as a MC with transition matrix $P = D^{-1}A$ where $D = Diag(Ae) = Diag(d)$. We assume that every node has at least one out-going edge, which can include self-loops. Note that $D_{ii} = d_i$, the degree of node $i$. The graph is strongly connected implies that the associated MC is irreducible with $p_{ij} = 1/d_j$ for those states $j$ such that $i \to j$, and 0 otherwise. The graph is undirected implies that the associated MC is reversible, and that the stationary probability vector is $\pi^T = d/d^T e$.

Kirkland (2010) considers mixing on directed graphs. For any stochastic matrix $P$ of order $m$, the directed graph associated with $P$, $D(P)$ is the directed graph on the vertices labelled 1, 2, ..., $m$ such that for each $i,j$ = 1, 2, ..., $m$, $i \to j$ is an arc on $D(P)$ if and only if $p_{ij} > 0$. For a strongly connected graph $D$ on $m$ vertices define the class $\sum_D = \{P | P$ is stochastic and $m \times m$ and for each $i, j$ = 1, 2, ..., $m$, $i \to j$ is an arc on $D(P)$ only if $i \to j$ is an arc in $D\}$. Define Kemeny's constant $K(P)$ with the convention that $m_{ii} = 0$. Let $\mu(D) = \inf \{K(P) | P \in \sum_D$ and $P$ has 1 as a simple eigenvalue$\}$. Let $k$ = the length of the longest cycle in $D$, then $\mu(D) = (2m - k - 1)/2$. Thus, if the MC periodic with period d = $m$, $\mu(D) = (m-1)/2$.

**Electric networks and graphs**
Doyle and Snell, (1984) in a very readable text established that there is a connection between electric networks and random walks (RWs) and graphs. On a connected graph $G$ with vertex set $V$ = $\{1,2, ..., m\}$ assign to the edge $(i, j)$ a resistance $r_{ij}$. The conductance of an edge $(i, j)$ is $C_{ij} = 1/r_{ij}$. Define a RW on $G$ to be a MC with transition probabilities $p_{ij} = C_{ij}/C_i$ with $C_i = \sum_j C_{ij}$. The graph is connected implies that the MC is ergodic with a stationary probability vector $\pi^T = (\pi_1,...,\pi_m)$ where $\pi_j = C_j/C$ with $C = \sum_i C_i$. The MC is in fact reversible. On the electric network we define $C_{ij} = \pi_i p_{ij}$. (If $p_{ii} \neq 0$ the resulting network will need a conductance from $i$ to $i$.)

For a network of resistors assigned to the edges of a connected graph choose two points $a$ and $b$ and put a one volt battery across these points establishing a voltage $v_a = 1$, $v_b = 0$. We are interested in finding the voltages $v_i$ and the currents $I_{ij}$ in the circuit and to give a probabilistic interpretation. By Ohm's Law, $I_{ij} = (v_i - v_j)/r_{ij}$ = $(v_i - v_j)C_{ij}$. Note $I_{ij} = -I_{ji}$. By Kirchhoff's current law $\sum_j I_{ij} = 0$ for $i \neq a, b$. i.e. if $\sum_j (v_i - v_j)C_{ij} = 0$ implying $v_i = \sum_j v_j p_{ij}$ for $i \neq a, b$. Let $h_i$ be the probability of starting at $i$, that state $a$ is reached before $b$. Then $h_i$ also satisfies above equations with $v_a = h_a = 1$ and $v_b = h_b = 0$. i.e. we can interpret the voltage as a "hitting probability".

Let $E_a T_b$ be the expected value, starting from the vertex $a$, of the hitting time $T_b$ of the vertex $b$. Let $\pi_i$ be the stationary probability of the MC at vertex $i$. When we impose a voltage $v$ between points $a$ and $b$ a voltage $v_a = v$ is established at $a$ and $v_b = 0$ and a current $I_a = \sum_x I_{ax}$ will flow into the circuit from outside the source.

We define the effective resistance between $a$ and $b$ as $R_{ab} = v_a/i_a$, as calculated using Ohm's Law. Using these concepts, Palacios and Tetali, (1996) showed that

$$E_a T_b = \frac{1}{2} \sum_i C_i \{R_{ab} + R_{bi} - R_{ai}\}.$$

Let $G$ be a simple connected graph with vertex set $V = \{1, 2, ..., m\}$. Let $R_{ij}$ be the effective resistance between $i$ and $j$. The *Kirchhoff index* was defined by Klein and Randic (1993) as $Kf(G) = \sum_{i<j} R_{ij}$. Since $R_{ij} = R_{ji}$ and $R_{ii} = 0$, $Kf(G) = \frac{1}{2} \sum_{i,j} R_{ij}$. This index has been considered by a number of authors and has been used in Chemistry to discriminate between different molecules with similar shapes and cycle structures.

Let us use the relations between electric networks and RWs on graphs. For a graph of $m$ vertices computing $Kf(G)$ entails finding $O(m^2)$ values of the $R_{ij}$ that is equivalent to finding $O(m^2)$ values of the $E_i T_j$ for the RW on the graph. Palacois and Renom, (2010) showed that $Kf(G)$ can be characterised as

$$Kf(G) = \frac{1}{2|E|} \sum_{i,j} E_i T_j.$$

The argument that they used was based on the fact that the "commute times" can be expressed as $E_i T_j + E_j T_i = 2|E|R_{ij}$, as derived by Aldous and Fill, (2002).

Zhu, Klein, Lukovits, (1996); Gutman, Mohar, (1996); and Mohar, Babic, Trinajstic, (1993) showed that $Kf(G)$ can also be characterised as $Kf(G) = m \sum_{i=1}^{m-1} \frac{1}{\mu_i}$ (where the $\mu_i$'s ($i = 1, 2,.., m$) with $\mu_m = 0$ are the eigenvalues of the (ordinary or combinatorial) Laplacian matrix $L$ of $G$, i.e. $L = D - A = D(I - P)$.

Using the above characterisation, Zhou and Trinajstic, (2009), found upper and lower bounds for $Kf(G)$ in terms of the eigenvalues of the normalised Laplacian $L = D^{-1/2} L D^{-1/2}$.

We have mentioned these results because of the recent research of Palacois (2010) establishing a connection between the Kirchhoff index and $Z$, the fundamental matrix of the underlying MC. In the case of $d$-regular graphs, (where all vertices have exactly $d$ neighbours) and using the characterisation of the Kirchhoff index as $Kf(G) = \frac{1}{d} \sum_j E_1 T_j$, Palacois showed that $Kf(G) = \frac{m}{d} [tr(Z) - 1]$ where $Z = (I - P + e\pi^T)^{-1}$, with $P$ the transition matrix of the random walk and $\pi^T$ its stationary probability vector.

Thus we have a connection between the Kirchhoff index and the modified Kemeny's constant $K = tr(Z) - 1$. This opens up a new line of research in an attempt to establish alternative interpretations to Kemeny's constant. There has been a lot of interest in recent years in graph theory, Laplacian and normalised Laplacians, electric networks, and hitting times and there is further scope to explore the interconnections alluded to above.

**Variances of mixing times**

We conclude this survey by considering a possible extension. We saw that the expected time to "mixing" starting in any state is $K$, a constant independent of the starting state. What about the variance of the mixing times? Do these depend on the starting state? If so, can we choose a desirable starting state? In Hunter (2008) these questions were explored. The results depended on expressions for the second moments of the first passage time variables.

Let $m_{ij}^{(2)}$ be the 2-*nd* moment of the first passage time from state *i* to state *j*. i.e. $m_{ij}^{(2)} = E[T_{ij}^2 \mid X_0 = i]$ for all $(i, j)$ in S×S; and let $M^{(2)} = [m_{ij}^{(2)}]$.

Let *T* be the mixing time variable and let $\eta_i^{(k)} = E[T^k \mid X_0 = i] = \sum_{j=1}^{m} m_{ij}^{(k)} \pi_j$ and $\boldsymbol{\eta}^{(k)T} = (\eta_1^{(k)}, \eta_2^{(k)}, \ldots, \eta_m^{(k)})$.

We have seen that $\boldsymbol{\eta}^{(1)T} = K\boldsymbol{e}$, i.e. the expected mixing times, starting at *i*, are constant. The variance of the mixing time, starting at *i* is given by $v_i = \eta_i^{(2)} - K^2$. If $\boldsymbol{v}^T = (v_1, v_2, \ldots, v_m)$ then $\boldsymbol{v} = \boldsymbol{\eta}^{(2)} - K^2 \boldsymbol{e}$.

From Hunter (2008), if *G* is any g-inverse of $I - P$, such that $G\boldsymbol{e} = \boldsymbol{e}$,
$\boldsymbol{\eta}^{(2)} = [2tr(G^2) - 3tr(G) - (1 - 2g)(1 - g)]\boldsymbol{e} + 2L\boldsymbol{\alpha}$,
$\boldsymbol{v} = [2tr(G^2) - (tr(G))^2 - (5 - 2g)tr(G) - (1 - g)(2 - 3g)]\boldsymbol{e} + 2L\boldsymbol{\alpha}$,
where $L = I - G + EG_d$ and $\boldsymbol{\alpha} = \boldsymbol{e} - (\Pi G)_d D\boldsymbol{e} + G_d D\boldsymbol{e}$.
$v_i = v$ for all *i* if and only if $L\boldsymbol{\alpha} = l\boldsymbol{e}$. A sufficient condition is $\boldsymbol{\alpha} = \alpha\boldsymbol{e}$.

The special case of 2-states was considered in detail. Let $P = \begin{bmatrix} 1-a & a \\ b & 1-b \end{bmatrix}$ and $d = 1 - a - b$. Then it follows that

$$\boldsymbol{v} = \begin{bmatrix} v_1 \\ v_2 \end{bmatrix} = \frac{1}{ab(1-d)^2} \begin{bmatrix} (2a^2 + 2b - 3ab)(a+b) - ab \\ (2b^2 + 2a - 3ab)(a+b) - ab \end{bmatrix}$$

The equality of the variances occurs $v_1 = v_2$ if and only if either $a = b$ (symmetry) or $a = 1 - b$ (independent trials). Thus, except in these two situations, it is impossible for the mixing time variables to have identical variances.